\documentclass[11pt]{article}
\usepackage{amssymb}
\usepackage{mathrsfs}
\usepackage{latexsym,amsmath,amssymb,amsfonts,epsfig,graphicx,cite,psfrag}
\usepackage{eepic,color,colordvi,amscd}
\usepackage{ebezier}
\usepackage{verbatim}
\usepackage{subfigure}
\usepackage{enumitem}
\usepackage{booktabs}
\usepackage{amsthm}
\usepackage{comment}
\usepackage{color}
\theoremstyle{plain}
\newtheorem{theo}{Theorem}[section]

\newtheorem{lem}[theo]{Lemma}
\newtheorem{coro}[theo]{Corollary}

\newtheorem{conj}[theo]{Conjecture}

\theoremstyle{definition}
\newtheorem{defi}[theo]{Definition}

\theoremstyle{remark}

\def\qed{\hfill \rule{4pt}{7pt}}

\newcommand{\OPT}{\operatorname{OPT}}

\newcommand{\re}{\operatorname{Re~}}

\def\mi{\textbf{i}}

\let\svthefootnote\thefootnote
\newcommand\blankfootnote[1]{%
	\let\thefootnote\relax\footnotetext{#1}%
	\let\thefootnote\svthefootnote%
}

\begin{document}
\title{Minimizing cycles in tournaments and normalized $q$-norms}
	
\author{Jie Ma\thanks{School of Mathematical Sciences, University of Science and Technology of China, Hefei, Anhui, China.
Partially supported by NSFC grant 11622110, the project ``Analysis and Geometry on Bundles'' of Ministry of Science and Technology of the People's Republic of China,
and Anhui Initiative in Quantum Information Technologies grant AHY150200. Email: jiema@ustc.edu.cn.}~~~~~~~
Tianyun Tang\thanks{Department of Mathematics, National University of Singapore, Singapore. Email: ttang@u.nus.edu.}
}	

\date{}
	
\maketitle

\begin{abstract}
Akin to the Erd\H{o}s-Rademacher problem, Linial and Morgenstern \cite{C4} made the following conjecture in tournaments:
for any $d\in (0,1]$, among all $n$-vertex tournaments with $d\binom{n}{3}$ many 3-cycles, the number of 4-cycles is asymptotically minimized by a special random blow-up of a transitive tournament.
Recently, Chan, Grzesik, Kr\'al' and Noel \cite{C3C4} introduced spectrum analysis of adjacency matrices of tournaments in this study, and confirmed this for $d\geq 1/36$.

In this paper, we investigate the analogous problem of minimizing the number of cycles of a given length.
We prove that for integers $\ell\not\equiv 2\mod 4$, there exists some constant $c_\ell>0$ such that if $d\geq 1-c_\ell$,
then the number of $\ell$-cycles is also asymptotically minimized by the same family of extremal examples for $4$-cycles.
In doing so, we answer a question of Linial and Morgenstern \cite{C4} about minimizing the $q$-norm of a probabilistic vector with given $p$-norm for any integers $q>p>1$.
For integers $\ell\equiv 2\mod 4$, however the same phenomena do not hold for $\ell$-cycles,
for which we can construct an explicit family of tournaments containing fewer $\ell$-cycles for any given number of $3$-cycles.
We conclude by proposing two conjectures on the minimization problem for general cycles in tournaments.
\end{abstract}

\section{Introduction}
We consider extremal problems on cycles in tournaments.\footnote{Throughout this paper, by a cycle in a digraph, we always mean a {\it directed} one.}
For a cycle $C_{\ell}$ of length $\ell$ and a tournament $T$,
a {\it homomorphism} of $C_{\ell}$ to $T$ is a mapping from $V(C_{\ell})$ to $V(T)$ that preserves edges.
Let $t(C_{\ell},T)$ denote the number of homomorphisms of $C_{\ell}$ to $T$ divided by $|V(T)|^{\ell}$.

One of the most natural extremal problems on tournaments is to determine the maximum number of cycles of a given length in tournaments.
This can be traced back to the work of Kendall and Babington Smith \cite{KS} in 1940.
It is well-known (see \cite{quasitournament,G59}) that a tournament $T$ satisfies $t(C_3,T)\leq \frac18$ and it has the maximum number of cycles $C_3$ if and only if it is almost regular.
Optimal results on 4-cycles and 5-cycles were obtained by Beineke and Harary \cite{BH} and by Komarov and Mackey \cite{KM}, respectively.
For other cycles, Day conjectured in \cite{Day} that the asymptotic maximum of $t(C_{\ell},T)$ is achieved by a random oriented tournament if and only of $\ell$ is not divisible by four.
This was confirmed recently by Grzesik, Kr\'al', Lov\'asz and Volec in \cite{Ckgood} using algebraic approach,
where they also obtained the asymptotic maximum of $t(C_8,T)$ and a very close estimation on the maximum of $t(C_{\ell},T)$ for any $\ell$ divisible by four.

Another natural extremal problem is to study the complementary problem, i.e., the minimum number of cycles of a given length in tournaments.
Since irregular tournaments (such as {\it transitive} tournaments) can contain arbitrary small number of cycles,
one has to introduce some extra restriction on tournaments for measuring regularity in this minimizing problem.
As mentioned the results of \cite{quasitournament,G59} earlier, one such good measure could be the density $t(C_3,T)$ of 3-cycles.
This was indeed addressed by Linial and Morgenstern in \cite{C4},
where they asked for the asymptotic minimum density $t(C_4,T)$ of 4-cycles in tournaments $T$ with fixed density $t(C_3,T)$ of 3-cycles.

A {\it random blow-up} of a $m$-vertex transitive tournament is a tournament $T$ with $V(T)=V_1\cup V_2\cup ...\cup V_m$
such that all arcs within each $V_i$ are oriented randomly and for any $i<j$, all arcs between $V_i$ and $V_j$ are oriented from $V_i$ to $V_j$.
Linial and Morgenstern \cite{C4} conjectured that a random blow-up $T^*$ of a transitive tournament with all but one part of equal size and one smaller part would achieve the asymptotic minimum $t(C_4,T)$.
Suppose that $T^*$ has $n$ vertices and $t$ parts of equal size $zn$ for some real $z\in (0,1]$,
then $tz\leq 1< (t+1)z$, implying that $t=\lfloor z^{-1}\rfloor$.
Therefore, with high probability, it holds for any $\ell\geq 3$ that
$$t(C_\ell, T^*)=\frac{1}{2^\ell}\Big(\lfloor z^{-1}\rfloor z^\ell+(1-\lfloor z^{-1}\rfloor z)^\ell\Big)+o(1), \mbox{ ~~ where } o(1)\to 0 \mbox { as } n\to \infty.$$
The authors \cite{C4} also pointed out that the structure of this conjectured extremal configuration resembles those of the famous Erd\H{o}s-Rademacher problem \cite{Erd}
on the minimum number of triangles in a graph of a given number of vertices and edges (see \cite{C3C4} for more discussion).

To state the above conjecture in a precise formula, define a function $g_\ell:[0,1/8]\rightarrow[0,1]$ for any integer $\ell\geq 4$ as follows:
Let $g_\ell(0)=0$ and for any real $z\in (0,1],$ let
\begin{align*}
&g_\ell\left(\frac18\Big(\lfloor z^{-1}\rfloor z^3+(1-\lfloor z^{-1}\rfloor z)^3\Big)\right)= \frac{1}{2^\ell}\Big(\lfloor z^{-1}\rfloor z^\ell+(1-\lfloor z^{-1}\rfloor z)^\ell\Big).
\end{align*}
It is worth noting that the function $g_\ell(\cdot)$ is continuous and increasing on $[0,1/8]$.
\begin{conj}[Linial and Morgenstern \cite{C4}, Conjecture 2.2]\label{conj1}
Every tournament $T$ satisfies that
$$t(C_4,T)\geq g_4(t(C_3,T))+o(1),$$
where the $o(1)$ term goes to zero as $|V(T)|$ goes to infinity.
\end{conj}

Linial and Morgenstern \cite{C4} proved several other general bounds between the densities of $C_3$ and $C_4$ in tournaments.
In particular, some of these imply that this conjecture holds for tournaments $T$ with $t(C_3,T)$ asymptotically equal to $0, 1/8$ and $1/32$.
They \cite{C4} also asked to understand the relationships among the higher densities $t(C_\ell, T)$.

Very recently, Chan, Grzesik, Kr\'al' and Noel \cite{C3C4} introduced spectral analysis of adjacency matrices of tournaments in this study.
They proved the following result, which makes a breakthrough towards Conjecture~\ref{conj1}.

\begin{theo}[Chan, Grzesik, Kr\'al' and Noel \cite{C3C4}]\label{C3C4theo}
Conjecture~\ref{conj1} holds for $t(C_3,T)\in [\frac{1}{72},\frac{1}{8}]$.
\end{theo}

Moreover, they \cite{C3C4} developed a limit theory for tournaments and used it to
characterize the asymptotic structure of all extremal tournaments $T$ with $t(C_3,T)\in [1/32, 1/8]$ for Conjecture~\ref{conj1}.

In this paper, we study the minimum number of cycles in tournaments $T$ with fixed density $t(C_3,T)$ of 3-cycles.
The following is our main result, which exhibits an analog of Theorem~\ref{C3C4theo}
for a cycle of length not of the form $4k+2$ in tournaments that are ``close'' to be regular.

\begin{theo}\label{4k+3}
Let $T$ be a tournament and $\ell\geq 4$ be an integer with $\ell\not\equiv 2\mod 4$.
If $t(C_3,T)\geq \frac18-\frac{1}{10\ell^2}$,
then it holds that
$$t(C_\ell,T)\geq g_\ell(t(C_3,T))+o(1),$$
where the $o(1)$ term goes to zero as $|V(T)|$ goes to infinity.
\end{theo}

By the same random blow-up example, we see that the above lower bound
for any cycle $C_\ell$ of length not of the form $4k+2$ is asymptotically tight.
We conjecture that similarly as in Conjecture~\ref{conj1}, the bound in Theorem~\ref{4k+3} can be extended to any value of $t(C_3,T)$.

The same minimizing problem for cycles of length of the form $4k+2$ is much more complicated.
For instance, it is still an open problem (see \cite{5and6}) to determine the asymptotic minimum density $t(C_6,T)$ of 6-cycles
even in {\it regular} tournaments $T$, i.e., $t(C_3,T)=1/8$.
On the other hand, we do know that the analogous bound in Theorem~\ref{4k+3} for these cycles does not hold for {\it any} value of $t(C_3,T)$ (see Lemma~\ref{lem:4k+2}).
For any $\ell$ of the form $4k+2$, we will describe an explicit family of tournaments $T$ with $t(C_\ell, T)< g_\ell(t(C_3,T))$ for any value of $t(C_3,T)$,
and conjecture that this family provides the optimal value for the minimum $t(C_\ell, T)$.
Using results of \cite{Ckgood}, one also can give a close estimation on the minimum $t(C_\ell, T)$ for regular tournaments,
which provides evidences for supporting the above conjecture; see Lemma~\ref{lem:regular-4k+2}.

We now turn back to Theorem~\ref{4k+3}.
The proof of Theorem~\ref{4k+3} is motivated by the one in \cite{C3C4}, which used spectral analysis on adjacency matrices of tournaments.
Through the spectral analysis, Theorem~\ref{4k+3} can be reduced to a list of optimization problems (see Subsection \ref{subsec:opt}),
whose constraints are non-linear non-convex polynomials of degree at most three and its objective function is a multinomial of degree $\ell$.
The base case of these optimization problems is a minimization problem on the interplay between $q$-norms of probability vectors (see Theorem~\ref{question} below).
A {\it probability vector} $\vec{w}$ is a vector consisting of entries $w_i$ for each integer $i\geq 1$, where all $w_i$'s are non-negative and add up to exactly one.
Linial and Morgenstern \cite{C4} proved that for any real $0<C<1$, the minimum of $\sum w_i^4$ among all probability vectors $\vec{w}$ satisfying $\sum w_i^3=C$
is attained by letting $w_1=...=w_{m}> w_{m+1}\geq 0$ and $w_i=0$ for any $i>m+1$,\footnote{Lemma 2.7 in \cite{C4} also requires that $m$ is the smallest integer among all possible such choices. However, as we shall explain in Section~\ref{sec:OPT} (in the proof of Theorem~\ref{questionfinite}) that $\vec{w}$ satisfying $\sum w_i=1, \sum w_i^p=C$, $w_1=\cdots=w_{m}> w_{m+1}\geq 0$ and $w_i=0$ for any $i\geq m+2$ is unique, which only depends on $p$ and $C$.}
which was used in the proof of Theorem~\ref{C3C4theo} in \cite{C3C4}.
Linial and Morgenstern \cite{C4} also raised the following ``natural sounding'' question:
{\it find the smallest $q$-norm among all probability vectors of given $p$-norm for any integers $q>p>1$}.
We answer this in the following strengthening, which holds for any reals $q>p>1$.
(We should postpone the definition of the function $f_{p,q}$ in Section \ref{sec:OPT}.)
\begin{theo}\label{question}
For any reals $q>p>1$ and $C\in (0,1)$, consider all probability vectors $\vec{w}$ satisfying $\sum w_i^p=C$.
Then the minimum of $\sum w_i^q$ among all such vectors equals $f_{p,q}(C)$,
which is attained by letting $w_1=w_2=\cdots=w_{m}> w_{m+1}\geq 0$ and $w_i=0$ for any $i\geq m+2$.\footnote{Perhaps it is worth pointing out that it is not trivial to see that there always exists an optimal vector to this infinite-dimensional optimization problem.}
\end{theo}

After settling this base case in the list of optimization problems of Theorem~\ref{4k+3},
we apply mathematical induction to solve other optimization problems.
The difficulty we face in the inductive step is that unlike in \cite{C3C4}, we cannot use the method of Lagrange multipliers
because we find that for longer cycles, the Lagrange points of these optimization problems are not unique and their forms are very complicated due to the high degree multinomial objective function.
To overcome this, we focus on a shorter interval for $t(C_3,T)$ and turn to some analytic method to complete the proof of Theorem~\ref{4k+3}.

The rest of the paper is organized as follows.
In Section~\ref{sec:OPT}, we prove Theorem~\ref{question}.
In Section~\ref{sec:cycles}, we prove Theorem~\ref{4k+3}.
In Section~\ref{sec:4k+2}, we consider cycles $C_\ell$ for $\ell=4k+2$,
by giving a family of tournaments $T$ with $t(C_\ell, T)< g_\ell(t(C_3,T))$ for every value of $t(C_3,T)$ (see Lemma~\ref{lem:4k+2})
and an asymptotic form of the minimum $t(C_\ell, T)$ as $\ell\to \infty$ for regular tournaments $T$ (see Lemma~\ref{lem:regular-4k+2}).
We conclude with two conjectures concerning the minimization problem discussed here.

\section{The optimization problem of Linial-Morgenstern}\label{sec:OPT}
In this section we prove Theorem~\ref{question}, which considers an infinite-dimensional minimization problem.
We should first prove an analogous finite-dimensional version (see Theorem~\ref{questionfinite}), and then show how to use this finite-dimensional version to derive Theorem~\ref{question}.
We also give a direct corollary of Theorem~\ref{question} at the end of this section.

Let $q>p>1$ be two reals. We define the function $f_{p,q}(\cdot): [0,1]\to [0,1]$ as follows:
Let $f_{p,q}(0)=0$ and for any $z\in (0,1]$, let
$$f_{p,q}\left(\lfloor z^{-1}\rfloor z^p+(1-\lfloor z^{-1}\rfloor z)^p\right)=\lfloor z^{-1}\rfloor z^q+(1-\lfloor z^{-1}\rfloor z)^q.$$
We point out that for any real $p>1$, $z\mapsto \lfloor z^{-1}\rfloor z^p+(1-\lfloor z^{-1}\rfloor z)^p$ is a strictly increasing continuous bijection from $(0,1]$ to $(0,1]$.
From this, it is easy to see that $f_{p,q}(\cdot)$ is a strictly increasing continuous function on $[0,1]$.
Also let us notice that $g_\ell(s)=f_{3,\ell}(8s)/2^\ell$ for any $s\in [0,1/8]$.

\subsection{A finite-dimensional version}\label{subsec:finite}

\begin{theo}\label{questionfinite}
For any reals $q>p>1$, $C\in (0,1)$ and positive integer $k$ with $\frac{1}{k^{p-1}}\leq C$,
consider all vectors $\vec{w}=(w_1,...,w_k)$ satisfying $\sum_{i=1}^{k}w_i=1, \sum_{i=1}^{k}w_i^p=C$ and $w_i\geq 0$ for each $i\in [k]$.
Then the minimum of $\sum_{i=1}^{k}w_i^q$ among all such vectors equals $f_{p,q}(C)$,
which is attained by letting $w_1=w_2=\cdots=w_{m}> w_{m+1}\geq 0$ and $w_i=0$ for any $i\geq m+2$.
\end{theo}

First, let us note that optimal solutions to the minimization problem in Theorem~\ref{questionfinite} always exist, by the following two facts:
\begin{itemize}
\item [(i)] Since the function $(w_1,\dots,w_k)\rightarrow\sum_{i=1}^{k}w_i^p$ is continuous, and $\frac{k}{k^p}\leq C< 1$,
there exists a vector $(w_1,w_2,\cdots,w_k)$ satisfying $\sum_{i=1}^{k}w_i=1$, $\sum_{i=1}^{k}w_i^p=C$ and $w_i\geq 0$ for each $i\in [k]$.
Thus, the feasible region is not empty.
\item [(ii)] It is clear that the feasible region for $\vec{w}$ is compact.
\end{itemize}
	
To prove Theorem~\ref{questionfinite}, we need to establish some lemmas.
The first one is a generalization of the well-known fact that the determinant of a Vandermonde matrix is non-zero.
\begin{lem}\label{matrix}
For reals $0<c_1<c_2<...<c_n$ and $s_1<s_2<...<s_n$, let $a_{i,j}=c_j^{s_i}$ for any $1\leq i, j\leq n$.
Then the rank of the matrix $\{a_{i,j}\}_{n\times n}$ is $n$.
\end{lem}
\begin{proof}
We prove it by induction on $n$. First, it holds trivially when $n=1$.	
Now suppose it holds for $n\leq k-1$ for some $k\geq 2$
and assume by contradiction that when $n=k$, there exists such a matrix $A=\{c_j^{s_i}\}_{k\times k}$ whose rank is less than $k$.
Then there exists a non-zero column vector $\vec{\lambda}=(\lambda_1,...,\lambda_k)^T$ such that $\lambda^T A=(0,0,...,0)$.
Thus, the function $\sum_{i=1}^{k}\lambda_ix^{s_i}$ has $k$ distinct zero points $c_1, c_2,..., c_k$.
This also means that $\sum_{i=1}^{k}\lambda_ix^{s_i-s_1}$ has $k$ distinct zero points $c_1, c_2,..., c_k$.
By Lagrange's mean value Theorem, the function $\sum_{i=2}^{k}\lambda_i(s_i-s_1)x^{s_i-s_1-1}$ has $k-1$ distinct zero points $d_1,..., d_{k-1}$, where $d_i\in (c_i,c_{i+1})$.
Let $\mu_i=\lambda_{i+1}(s_{i+1}-s_1)$ and $t_i=s_{i+1}-s_1-1$ for each $i\in [k-1]$.
Then for the matrix $B=\{d_j^{t_i} \}_{(k-1)\times (k-1)}$ and the vector $\vec{\mu}=(\mu_1,...,\mu_{k-1})^T$, we have $\vec{\mu}^T B=(0,0,...,0)$.
If $\vec{\mu}=\vec{0}$, then we see $\lambda_i=0$ for any $i\in \{2,3,...k \}$.
As $\sum_{i=1}^{k}\lambda_ix^{s_i}$ has $k$ distinct zero points, it forces $\lambda_1=0$ and thus $\vec{\lambda}=\vec{0}$, which is a contradiction.
Therefore $\vec{\mu}\neq \vec{0}$, which implies that the rank of $B$ is less than $k-1$.
This is contradictory to our induction hypothesis.
\end{proof}
	
Let $\vec{w}=(w_1,...,w_k)$ denote an optimal vector to the minimization problem in Theorem~\ref{questionfinite} such that $w_1\geq w_2\geq \cdots \geq w_k\geq 0$.
The next lemma says that there are at most two distinct positive real numbers that the entries of $\vec{w}$ can take.
\begin{lem}\label{<3}
There exist two reals $1>a>b>0$ such that for any $i\in [k]$, $w_i\in \{0,b,a\}$
\end{lem}
	
\begin{proof}
Without loss of generality, we may assume by contradiction that there are $w_1>w_2>w_3>0$.
The method of Lagrange multipliers says that, at the optimal vector, the gradient of the objective function is perpendicular to the tangent plane of the feasible region.
This implies that the rank of the following matrix $$ \begin{pmatrix}1&1&1\\w_1^{p-1}&w_2^{p-1}&w_3^{p-1}\\w_1^{q-1}&w_2^{q-1}&w_3^{q-1}\end{pmatrix}$$ is less than 3.
By Lemma~\ref{matrix}, we get a contradiction, completing this proof.
\end{proof}

The following lemma says more about the optimal vector $\vec{w}$ and is key for Theorem~\ref{questionfinite}.
	
\begin{lem}\label{abb}
There exist no reals $1>a>b>0$ and integers $i,j,\ell\in [k]$ such that $w_i=a,w_j=b$ and $w_\ell=b$.
\end{lem}

\begin{proof}
We assume by contradiction that $w_i=a, w_j=b$ and $w_\ell=b$.
Our proof strategy is to find a new vector in the feasible region which reduces the objective function.
We consider two cases separately, depending on whether $2(\frac{a+2b}{2})^p\geq a^p+2b^p$ or not.

\medskip

\noindent{\bf Case 1.} $2(\frac{a+2b}{2})^p\geq a^p+2b^p$.

\medskip	

Let $f:[0,1)\rightarrow \mathbf{R}$ be such that $f(x)=(a-2x)^p+2(b+x)^p$. We have $f'(x)=2p(b+x)^{p-1}-2p(a-2x)^{p-1}$.
So $f'(x)<0$ on $[0,\frac{a-b}{3})$ and $f'(x)>0$ on $(\frac{a-b}{3},\frac{a}{2}]$.
We also have $f(0)=a^p+2b^p$ and $f(\frac{a}{2})=2(\frac{a+2b}{2})^p\geq a^p+2b^p$.
Therefore, there exists one and only one $x_0\in (\frac{a-b}{3},\frac{a}{2}]$ such that $f(x_0)=a^p+2b^p$.
Let $c=a-2x_0$, $d=b+x_0$. Then we have
\begin{align}\label{equ:cd}
0\leq c<d, ~ c+2d=2b+a \mbox{ ~and~ } c^p+2d^p=a^p+2b^p.
\end{align}
		
We claim that $0\leq c<b<d<a$.
To see this, first note that by the strict convexity of $x\mapsto x^p$, we have $f(\frac{a-b}{2})=b^p+2(\frac{a+b}{2})^p<b^p+(b^p+a^p)=f(0)$.
Then $f(x)<a^p+2b^p=f(x_0)$ on $(0,\frac{a-b}{2}]$. This implies $x_0>\frac{a-b}{2}$ and thus $c<b<d$.
It remains to show $d<a$.
If $\frac{a+2b}{2}<a$, then $d=b+x_0\leq b+\frac{a}{2}<a$, as desired.
So we may assume $\frac{a+2b}{2}\geq a$, i.e., $a-b\leq \frac{a}{2}$.
Using the strict convexity of $x\mapsto x^p$ again, $f(a-b)=(2b-a)^p+2a^p> 2b^p+a^p=f(x_0)$.
Because $f'(x)>0$ on $[a-b,\frac{a}{2}]$, we have $f(x)\geq f(a-b)>f(x_0)$ for all $x\in [a-b,\frac{a}{2}]$.
This implies that $x_0<a-b$ and thus $d=b+x_0<a$, proving the claim.

We want to replace $(a,b,b)$ with $(d,d,c)$ in $\vec{w}$ while reducing the objective function.
In order to do this, we consider the following function $h(x)$:
$$\mbox{Let } h(x)=c^x+2d^x-a^x-2b^x \mbox{ if } c>0, \mbox{ and let } h(x)=2d^x-a^x-2b^x \mbox{ if } c=0.$$

Now we argue that $h(q)<0$.
Assume by contradiction that $h(q)\geq 0$.
Using the fact that $0\leq c<b<d<a$, when $x$ is sufficiently large, we have $h(x)<0$.
So there exists some real $s\geq q$ such that $h(s)=0$.
By \eqref{equ:cd}, we also have $h(p)=h(1)=0$.
If $c>0$, then $h(0)=0$ and thus $h(x)$ has four zero points $0<1<p<s$.
This implies that $$\begin{pmatrix}1&1&1&1\\c&b&d&a\\c^p&b^p&d^p&a^p\\c^s&b^s&d^s&a^s\end{pmatrix}\begin{pmatrix}1\\-2\\2\\-1\end{pmatrix}=\begin{pmatrix}0\\0\\0\\0\end{pmatrix}.$$
Then the rank of the matrix above is less than 4, which is contradictory to Lemma~\ref{matrix}.
So $c=0$. Then $h(x)$ has three zero points $1<p<s$, which imply that $$\begin{pmatrix}b&d&a\\b^p&d^p&a^p\\b^s&d^s&a^s\end{pmatrix}\begin{pmatrix}-2\\2\\-1\end{pmatrix}=\begin{pmatrix}0\\0\\0\end{pmatrix}.$$
Again, this is contradictory to Lemma~\ref{matrix}.
We now complete the proof of $h(q)<0$.

Therefore, $c^q+2d^q<a^q+2b^q$.
Let $\vec{w}'$ be obtained from $\vec{w}$ by replacing $(a,b,b)$ with $(c,d,d)$.
Clearly, $\vec{w}'$ is still in the feasible region, but it reduces the objective function.
This contradicts that $\vec{w}$ is an optimal vector, completing the proof for Case 1.
	
\medskip	

\noindent{\bf Case 2.} $2(\frac{a+2b}{2})^p<a^p+2b^p$.

\medskip	

We first prove that $a>2b$.
Let $y(x)=x^p+2-2(\frac{x}{2}+1)^p$. So $y'(x)=px^{p-1}-p(\frac{x}{2}+1)^{p-1}$.
This shows that $y(x)$ is monotone decreasing on $(0,2)$ and monotone increasing on $(2,+\infty)$.
So $y(x)<y(0)=0$ for any $x\in (0,2]$.
Using the condition of Case 2, $y(\frac{a}{b})=(\frac{a}{b})^p+2-2(\frac{a}{2b}+1)^p>0$,
which implies that $\frac{a}{b}>2$. So $a>2b$ and moreover, $a+2b>a>\frac{a+2b}{2}$.
		
Let $g(x)=x^p+(a+2b-x)^p-a^p-2b^p$.
Since $g'(x)=px^{p-1}-p(a+2b-x)^{p-1}$, it is easy to see that $g(x)$ is monotone increasing on $[\frac{a+2b}{2},a+2b]$.
Because $g(a)=a^p+(2b)^p-a^p-2b^p>0$ and $g(\frac{a+2b}{2})=2(\frac{a+2b}{2})^p-a^p-2b^p<0$,
there exists a unique $x_0\in (\frac{a+2b}{2},a)$ such that $g(x_0)=0$.
Let $c=a+2b-x_0$ and $d=x_0$. Then one can derive that
\begin{align}\label{equ:cd2}
0<b<c<d<a, ~ c+d=2b+a \mbox{ ~and~ } c^p+d^p=a^p+2b^p.
\end{align}

In this case, we want to replace $(a,b,b)$ with $(d,c,0)$ in $\vec{w}$, while reducing the objective function.
For this purpose, we consider the function $h(x)=(\frac{c}{a})^x+(\frac{d}{a})^x-1-2(\frac{b}{a})^x.$

We claim that $h(x)>0$ holds for any $x\in (1,p)$.
Assume by contradiction that there exists $s\in (1,p)$ such that $h(s)\leq 0$.
By \eqref{equ:cd2}, $h(1)=h(p)=0$.
Let $\beta$ be a minimum point of $h(x)$ in $[1,p]$.
Because there exists $s\in (1,p)$ with $h(s)\leq 0$, we can choose $\beta\in (1,p)$ with $h(\beta)\leq 0$ and $h'(\beta)=0$.
Let $\alpha\in [0,\beta]$ be a maximum point of $h(x)$ on $[0,\beta]$.
As $h(0)=-1<0$ and there exists some $s'\in (0,\beta)$ with $h(s')\geq 0\geq \max\{h(0), h(\beta)\}$,\footnote{Clearly, here one can take $s'=1$.} we can choose $\alpha\in (0,\beta)$ with $h'(\alpha)=0$.
As $\lim_{x\rightarrow +\infty}h(x)=-1<0$, there exists $t>p$ such that $h(t)<0$.
Let $\gamma$ be a maximum point of $h(x)$ on $[\beta,t]$.
Since $h(p)\geq 0\geq \max\{h(\beta), h(t)\}$ where $p\in (\beta, t)$,
we can choose $\gamma\in (\beta,t)$ with $h'(\gamma)=0$.
For $0<\alpha<\beta<\gamma$, we have $h'(\alpha)=h'(\beta)=h'(\gamma)=0$, implying that
$$\begin{pmatrix}(\frac{b}{a})^{\alpha}&(\frac{c}{a})^{\alpha}&(\frac{d}{a})^{\alpha}\\
(\frac{b}{a})^{\beta}&(\frac{c}{a})^{\beta}&(\frac{d}{a})^{\beta}\\
(\frac{b}{a})^{\gamma}&(\frac{c}{a})^{\gamma}&(\frac{d}{a})^{\gamma}\end{pmatrix}\begin{pmatrix}-2\log\frac{b}{a}\\
\log\frac{c}{a}\\\log\frac{d}{a}\end{pmatrix}=\begin{pmatrix}0\\0\\0\end{pmatrix}.$$
Then the rank of the matrix above is less than 3, a contradiction to Lemma~\ref{matrix}.

Next we show that $h(x)<0$ holds for any $x\in (p,+\infty)$.
Assume by contradiction that there exists some $s\in (p,+\infty)$ such that $h(s)\geq 0$.
Because $\lim_{x\rightarrow +\infty}h(x)=-1$, there exists $t\in (s,+\infty)$ such that $h(t)<0$.
Let $\alpha$ be a maximum point of $h(x)$ on $[1,p]$.
By the previous paragraph, we can choose $\alpha\in (1,p)$ with $h(\alpha)>0$ and $h'(\alpha)=0$.
Let $\gamma$ be a maximum point of $h(x)$ on $[p,t]$.
Because $h(s)\geq 0\geq \max\{h(p),h(t)\}$ where $s\in (p,t)$, we can choose $\gamma\in (p,t)$ with $h(\gamma)\geq 0$ and $h'(\gamma)=0$.
Lastly, let $\beta$ be a minimum point of $h(x)$ on $[\alpha,\gamma]$.
Similarly, we can choose $\beta\in (\alpha,\gamma)$ with $h'(\beta)=0$.
For $0<\alpha<\beta<\gamma$, we have $h'(\alpha)=h'(\beta)=h'(\gamma)=0$.
This implies that $$\begin{pmatrix}(\frac{b}{a})^{\alpha}&(\frac{c}{a})^{\alpha}&(\frac{d}{a})^{\alpha}\\ (\frac{b}{a})^{\beta}&(\frac{c}{a})^{\beta}&(\frac{d}{a})^{\beta}\\(\frac{b}{a})^{\gamma}&(\frac{c}{a})^{\gamma}&(\frac{d}{a})^{\gamma}\end{pmatrix}\begin{pmatrix}-2\log\frac{b}{a}\\\log\frac{c}{a}\\\log\frac{d}{a}\end{pmatrix}=\begin{pmatrix}0\\0\\0\end{pmatrix}.$$
By Lemma~\ref{matrix} again, we get a contradiction.

In particular, $h(q)<0$, which means $c^q+d^q<a^q+2b^q$.
Then we can replace $(a,b,b)$ with $(d,c,0)$ in $\vec{w}$.
The new vector is still in the feasible region but reduces the objective function.
Thus we get a contradiction. This completes the proof of Lemma~\ref{abb}.	
\end{proof}

Now we are ready to prove Theorem~\ref{questionfinite}.

\begin{proof}[Proof of Theorem~\ref{questionfinite}.]
Let $\vec{w}=(w_1,...,w_k)$ be an optimal vector to the minimization problem in Theorem~\ref{questionfinite} such that $\sum_{i=1}^k w_i=1, \sum_{i=1}^k w_i^p=C$ and $w_1\geq w_2\geq \cdots \geq w_k\geq 0$.

By Lemmas~\ref{<3} and \ref{abb}, we may assume that there exists some $m$ such that $w_1=\cdots=w_{m}> w_{m+1}\geq 0$ and $w_i=0$ for any $i\geq m+2$.
We show that $\vec{w}$ satisfying the above additional condition is uniquely determined by $p$ and $C$.
Set $z=w_1$. Since $0\leq w_{m+1}<z$, we have $mz\leq 1=\sum w_i<(m+1)z$, implying that $z^{-1}-1<m\leq z^{-1}$, i.e., $m=\lfloor z^{-1}\rfloor$.
So $w_{m+1}=1-\lfloor z^{-1}\rfloor z$. Using the restriction $\sum_i w_i^p=C$, we get that
\begin{align}\label{equ:z}
\lfloor z^{-1}\rfloor z^p+(1-\lfloor z^{-1}\rfloor z)^p=C.
\end{align}
Note that $z\mapsto \lfloor z^{-1}\rfloor z^p+(1-\lfloor z^{-1}\rfloor z)^p$ is a strictly increasing continuous bijection from $(0,1]$ to $(0,1]$.
So there exists a unique solution $z\in (0,1)$ to the above equation. This proves that such an optimal vector $\vec{w}$ is unique,
which gives the minimum of $\sum_{i=1}^k w_i^q$ to be $\lfloor z^{-1}\rfloor z^q+(1-\lfloor z^{-1}\rfloor z)^q=f_{p,q}(C)$, completing the proof of Theorem~\ref{questionfinite}.
\end{proof}

\subsection{Proof of Theorem~\ref{question}}
We move on to prove Theorem~\ref{question}.
The proof will use the fact that the minimum value $f_{p,q}(C)$ of Theorem~\ref{questionfinite} is an increasing continuous function with the variable $C$ on $[0,1]$.
	
\begin{proof}[Proof of Theorem~\ref{question}.]
Consider any probability vector $\vec{w}$ satisfying $\sum w_i^p=C$.
Take any small constant $\epsilon>0$ and choose $k\in {\bf N}^+$ to be sufficiently large such that $t:=\sum_{i=1}^{k}w_i\geq 1-\epsilon$ and $\frac{1}{k^{p-1}}\leq C/2\leq C-\epsilon$.
Then we deduce that $\sum_{i=k+1}^{\infty}w_i\leq \epsilon$ and $\sum_{i=k+1}^{\infty}w_i^p\leq \epsilon$, implying that $\sum_{i=1}^{k}w_i^p\geq C-\epsilon$.
Now we normalize $w_i$ by letting $v_i:=w_i/t$ for each $i\in [k]$.
Then $\sum_{i=1}^{k}v_i=1$ and as $t\leq 1$, $\sum_{i=1}^{k}v_i^p=\left(\sum_{i=1}^{k}w_i^p \right)/t^p\geq C-\epsilon$.
By Theorem~\ref{questionfinite}, we have
$$\sum_{i=1}^{\infty}w_i^q\geq \sum_{i=1}^{k}w_i^q=t^q\cdot \sum_{i=1}^{k}v_i^q\geq t^q\cdot f_{p,q}\left(\sum_{i=1}^{k}v_i^p\right)\geq (1-\epsilon)^q \cdot f_{p,q}\left(C-\epsilon \right).$$
Using the continuity of $f_{p,q}$, by letting $\epsilon\rightarrow 0$ we can get $\sum_{i=1}^{\infty}w_i^q\geq f_{p,q}(C)$.
Now let $\vec{w}_0$ be a probability vector with $z=w_1=\cdots=w_{m}> w_{m+1}\geq 0$ and $w_i=0$ for any $i\geq m+2$,
where $m=\lfloor z^{-1}\rfloor$ and $z$ is from \eqref{equ:z}.
Then such $\vec{w}$ satisfies that $\sum_{i=1}^{\infty}w_i^p=C$ and $\sum_{i=1}^{\infty}w_i^q= f_{p,q}(C)$.
So the minimum value in Theorem~\ref{question} exists and equals $f_{p,q}(C)$.
\end{proof}
	
Now we present a corollary of Theorem~\ref{question}, which will be used frequently in Section~\ref{sec:cycles}.

\begin{coro}\label{coro}
Let $n, \ell$ be positive integers and $x_1,\dots,x_n$ be any non-negative reals such that $\ell\geq 4$ and $\sum_{i=1}^n x_i=\frac{1}{2}$.
Then $\sum_{i=1}^{n} x_i^\ell\geq g_\ell\Big(\sum_{i=1}^{n}x_i^3\Big).$
\end{coro}
\begin{proof}
By definition of $g_\ell(\cdot)$, together with Theorem~\ref{question}, we have that for any $s\in [0,1/8]$,
\begin{align*}
g_\ell(s)&=\frac{f_{3,\ell}(8s)}{2^\ell}=\min\left\{\frac{1}{2^\ell}\sum_i w_i^{\ell}\ |\ \sum_i w_i^3=8s,\
\sum_i w_i=1,\ w_i\geq 0 \ \mbox{for any}\ i\geq 1 \ \right\} \\
&= \min\left\{\sum_i x_i^{\ell}\ |\ \sum_i x_i^3=s,\
\sum_i x_i=\frac12,\ x_i\geq 0 \ \mbox{for any}\ i\geq 1 \ \right\}.
\end{align*}
This implies the desired inequality.
\end{proof}

\section{Cycles of given length in tournaments}\label{sec:cycles}

In this section, we prove Theorem~\ref{4k+3} by using the approach based on the spectrum analysis on the adjacency matrix of a tournament, which was initialed in \cite{C3C4,Ckgood}.

\subsection{Some preliminaries}

First we give some algebraic notation on a tournament $T$. Let $V(T)=[n]$.
For any $i,j\in [n]$, we write $ij$ for an arc with head $i$ and tail $j$.
The {\it adjacency matrix} of $T$ is an $n\times n$ matrix $M=(m_{i,j})$,
where $m_{i,j}=1$ if $ij\in E(T)$ and $m_{i,j}=0$ otherwise.
It is easy to see that
$$t(C_\ell,T)=tr(M^\ell)/n^\ell,$$
where $tr(X)$ denotes the trace of the matrix $X$.
Let $I$ denote the $n\times n$ identity matrix.
We call $A:=(\frac{1}{2}I+M)/n$ the {\it tournament matrix} of $T$.
It follows that $A+A^T=J$, where $J$ is the $n\times n$ matrix with each entry being $\frac{1}{n}$.
It is not difficult to see that
\begin{align}\label{equ:trA}
t(C_\ell,T)=tr(M^\ell)/n^\ell=tr(A^\ell)+o(1),
\end{align}
where $\lim_{n\rightarrow +\infty}o(1)=0$.
The following fact will be crucial (see \cite{TR}, or Lemma~2 in \cite{C3C4}).

\begin{lem}\label{Re}
Every eigenvalue of a tournament matrix has nonnegetive real part.
\end{lem}

Besides Corollary~\ref{coro}, we also need to use the following optimization.
The proof of the case $p=3$ was given in \cite{C3C4} (see Lemma~9), which can be easily transformed to general cases.

\begin{lem}\label{n}
For integers $p\geq 2$, $n\geq 1$ and any real $1>t>0$ such that $nt\geq 1$,
consider all vectors $(w_1,w_2,...,w_n) $ satisfying $\sum_{i=1}^{n}w_i=1$ and $w_i\in [0,t]$ for every $i\in [n]$.
Then the maximum of $\sum_{i=1}^{n}w_i^p$ is attained by letting $w_1=w_2=...=w_{\lfloor t^{-1}\rfloor}=t$, $w_{\lfloor t^{-1}\rfloor+1}=1-t\lfloor t^{-1}\rfloor$ and $w_j=0$ for any $j\geq \lfloor t^{-1}\rfloor+2$.
\end{lem}

For a complex number $z$, we denote its real part by $\re z$ and its absolute value by $|z|$.
We write the imaginary unit by $\mi$.
The following inequality was explicitly given in \cite{Ckgood} (see the inequality $(9)$ therein).
For completeness, we present a proof here.

\begin{lem}\label{Re inequality}
For any odd integer $\ell>0$ and complex number $z$ with $\re z\geq 0$, it holds that $\re z^\ell\leq \ell |z|^{\ell-1} \re z$.
\end{lem}
\begin{proof}
Let $z=a+b\mi$ for reals $a, b$ with $a\geq 0$.
It suffices to consider $a>0$ and $b\geq 0$. Write $b/a=\tan(\alpha)$ for $\alpha\in [0,\frac{\pi}{2})$.
Then $\re z^\ell=|z|^{\ell}\cdot \cos(\ell \alpha)$.
Since $\ell$ is odd, it follows that $|\re z^\ell|=|z|^{\ell}\cdot |\sin\ell (\pi/2-\alpha)|\leq |z|^{\ell}\cdot \ell \cdot |\sin (\pi/2-\alpha)|=|z|^{\ell}\cdot\ell\cdot |\cos \alpha|=\ell|z|^{\ell-1}a$.
\end{proof}

\subsection{Proof of Theorem~\ref{4k+3}}\label{subsec:opt}
In view of Theorem~\ref{C3C4theo}, we may assume that $\ell\geq 5$ is an integer with $\ell \not\equiv 2 \mod 4$.
Let $\ell=4k+\mu$, where $\mu\in \{-1,0,1\}$ and let $T$ be any tournament satisfying $$t(C_3,T)\geq \frac18-\frac{1}{10\ell^2}.$$
Our goal is to show $t(C_\ell,T)\geq g_\ell(t(C_3,T))+o(1),$ where $\lim_{|V(T)|\rightarrow +\infty}o(1)=0$.

Denote $A$ by its tournament matrix and let $\sigma=tr(A^3)$.
By the Perron-Frobenius theorem, we may assume that the eigenvalues of $A$ are $\rho,r_1,\dots,r_t,a_1+b_1\mi,\dots,a_s+b_s\mi$,
where $\rho$ is the spectral radius of $A$ and $r_i, a_j, b_j$ are reals for each $i\in [t]$ and $j\in [s]$.
We now collect some properties on these eigenvalues.
First, using Lemma~\ref{Re} together with the fact that $\rho$ is the spectral radius,
we see that $r_i\in [0,\rho]$ for each $i\in [t]$ and $a_j\geq 0$ and $a_j^2+b_j^2\leq \rho^2$ for each $j\in [s]$.
Since each diagonal element of $A$ is $\frac{1}{2n}$,
we have $\rho+\sum_{i=1}^{t}r_i+\sum_{j=1}^{s}a_j=\re tr(A)=tr(A)=\frac{1}{2}$.
Similarly, we have $$\rho^3+\sum_{i=1}^{t}r_i^3+\sum_{j=1}^{s}(a_j^3-3a_jb_j^2)=\re tr(A^3)=tr(A^3)=\sigma,$$
$$\rho^2+\sum_{i=1}^{t}r_i^2+\sum_{j=1}^{s}(a_j^2-b_j^2)=\re tr(A^2)=tr(A^2)\geq 0,$$
where the second inequality implies that $\sum_{j=1}^{s}b_j^2\leq \rho^2+\sum_{i=1}^{t}r_i^2+\sum_{j=1}^{s}a_j^2\leq (\rho+\sum_{i=1}^{t}r_i+\sum_{j=1}^{s}a_j)^2=\frac{1}{4}.$
Since $A$ is a real matrix, if $a_j+b_j\mi$ is an eigenvalue of $A$, then so is $a_j-b_j\mi$.
Hence we have $2b_j^2\leq \frac{1}{4}$, implying that $b_j^2\leq \frac18$ for each $j\in [s]$.

\begin{center}
	\begin{tabular}{ll}
        	\bottomrule
        	{} & \quad\quad Optimization Problem $\OPT_{\ell}^{\sigma}(s,t,\rho)$ \\
	\bottomrule
        	Parameters :&\quad\quad real numbers $\sigma\in [0,\frac{1}{8}]$, $\rho\in [0,\frac{1}{2}]$\\
	{}&\quad\quad nonnegative integers $\ell,t,s$ with $t+s\geq 1$ and $\ell\geq 5$\\
	Variables :&\quad\quad real numbers $r_1,\ldots,r_t,a_1,b_1,\dots,a_s,b_s$\\
	Constraints :&\quad\quad $\rho+\sum\limits_{i=1}^{t}r_i+\sum\limits_{j=1}^{s}a_j=\frac{1}{2}$\\
	{}&\quad\quad $\rho^3+\sum\limits_{i=1}^{t}r_i^3+\sum\limits_{j=1}^{s}(a_j^3-3a_jb_j^2)=\sigma$\\
	{}&\quad\quad $0\leq r_1,\ldots,r_t\leq \rho$, $0\leq a_1,\ldots,a_s$\\
	{}&\quad\quad for any $j\in [s]$, $a_j^2+b_j^2\leq \rho^2$, $b_j^2\leq \frac{1}{8}$\\
	Objective :&\quad\quad $\Phi_{\ell}^{\sigma}(s,t,\rho):=\min\Big\{\rho^\ell+\sum\limits_{i=1}^{t}r_i^\ell+\sum\limits_{j=1}^{s}{\re}(a_j+b_j\mi)^\ell\Big\}$\\
        \bottomrule
      	\end{tabular}
\end{center}	

We define the optimization problem $\OPT_{\ell}^\sigma(s,t,\rho)$, for which the optimal solution $\Phi_{\ell}^{\sigma}(s,t,\rho)$ exists.
Since the eigenvalues of $A$ satisfy the constraints in $\OPT_{\ell}^\sigma(s,t,\rho)$,
using \eqref{equ:trA} we have
\begin{align}\label{equ:tCl}
t(C_\ell,T)=tr(A^\ell)+o(1)=\rho^\ell+\sum_{i=1}^{t}r_i^\ell+\sum_{j=1}^{s}{\re}(a_j+b_j\mi)^\ell+o(1)\geq \Phi_{\ell}^{\sigma}(s,t,\rho)+o(1),
\end{align}
where $\lim_{|V(T)|\rightarrow +\infty}o(1)=0$.
Also by \eqref{equ:trA}, we have that
$\sigma=tr(A^3)=t(C_3,T)+o(1)\geq \frac18-\frac{1}{10\ell^2}+o(1)\geq \big(\frac{1}{2}-\frac{1}{80k^2}\big)^3+\big(\frac{1}{80k^2}\big)^3.$

\medskip

We now devote the rest of this section to show the following statement:
\begin{itemize}[itemindent=58pt]
\item [\bf Statement ($\star$).] For any real $\rho\in [0,\frac12]$ and nonnegative integers $\ell, t, s$ with $\ell\geq 5$, $\ell \not\equiv 2 \mod 4$ and $t+s\geq 1$,
if $\sigma\geq (\frac{1}{2}-\frac{1}{80k^2})^3+(\frac{1}{80k^2})^3$, then $\Phi_{\ell}^{\sigma}(s,t,\rho)\geq g_{\ell}(\sigma)$.	
\end{itemize}
Before we process, we point out that to finish the proof of Theorem~\ref{4k+3},
it suffices to prove ($\star$).
Indeed, if ($\star$) holds, then using \eqref{equ:tCl} and the fact that $g_\ell(\cdot)$ is a continuous function,
we can deduce the desired inequality as follows:
\begin{align*}
t(C_\ell,T)-o(1)\geq \Phi_{\ell}^{\sigma}(s,t,\rho)\geq g_\ell(\sigma)=g_\ell\big(t(C_3,T)+o(1)\big)\geq g_\ell(t(C_3,T))+o(1).
\end{align*}

We prove the statement ($\star$) by induction on $s$.
If $s=0$, then applying Corollary~\ref{coro} directly, we can derive $\Phi_{\ell}^{\sigma}(s,t,\rho)\geq g_{\ell}(\sigma)$.
So we may assume that $s\geq 1$ and there exists some integer $m$ such that ($\star$) holds for any $s<m$.
Now we consider when $s=m$. Assume that $(r_1,\dots,r_t,a_1,b_1,\dots,a_m,b_m)$ is an optimal vector of $\OPT_{\ell}^\sigma(m,t,\rho)$.

Recall $\ell=4k+\mu$, where $\mu\in \{-1, 0, 1\}$.
Suppose that $a_m=0$.
If $\mu=0$, then $\re (a_m+b_m\mi)^{4k}=b_m^{4k}$ is minimized in the objective function when $b_m=0$.
If $\mu\in \{-1,1\}$, then ${\re}(a_m+b_m\mi)^{\ell}={\re}(b_m\mi)^{\ell}=0$,
so by setting $b_m=0$ the new vector still satisfies all the constraints and the objective value remains unchanged.
Therefore, when $a_m=0$, we can always set $b_m=0$ and use induction directly.
From now on, we may assume $a_m>0$.

Since $\rho+\sum_{i=1}^t r_i+\sum_{j=1}^m a_j=\frac12$ and $0\leq r_i, a_j\leq \rho$ for each $i\in [t],j\in [m]$,
we can deduce from Lemma~\ref{n} that
$$\Big(\frac{1}{2}-\frac{1}{80k^2}\Big)^3+\Big(\frac{1}{80k^2}\Big)^3\leq \sigma \leq \rho^3+\sum_{i=1}^{t}r_i^3+\sum_{j=1}^{m}a_j^3 \leq\left\lfloor\frac{1}{2\rho}\right\rfloor\rho^3+\Big(\frac{1}{2}-\rho\left\lfloor\frac{1}{2\rho}\right\rfloor \Big)^3.$$
This implies that
\begin{align}\label{equ:am}
\rho\geq 1/2-1/80k^2 \mbox{ and thus } a_j\leq 1/80k^2 \mbox{ for each } j\in [m].
\end{align}

Let $\sigma':=\sigma+3a_mb_m^2$. Shifting $\sigma$ to $\sigma'$ and viewing $a_m$ as a new variable $r_{t+1}$,
we can obtain a new optimization problem $\OPT_{\ell}^{\sigma'}(m-1,t+1,\rho)$.
Our proof idea in what follows is to compare the optimal values of the two optimization problems $\OPT_{\ell}^\sigma(m,t,\rho)$ and $\OPT_{\ell}^{\sigma'}(m-1,t+1,\rho)$
by using analytic arguments.
For this purpose, we introduce some definitions as below.
For $|x|\leq |b_m|$, let $$\sigma(x)=\sigma+3a_mx^2 \mbox{~~ and ~~} \lambda(x)=\frac13\sqrt{6\sigma(x)-\frac{3}{16}}.$$
By the first two constraints of $\OPT_{\ell}^\sigma(m,t,\rho)$,
we have $$\frac{1}{32}\leq \sigma\leq \sigma(x)\leq \sigma(b_m)\leq\rho^3+\sum_{i=1}^{t}r_i^3+\sum_{j=1}^{m}a_j^3\leq \left(\rho+\sum_{i=1}^{t}r_i+\sum_{j=1}^{m}a_j\right)^3=\frac{1}{8}.$$
Solving $\left\lfloor z^{-1}\right\rfloor z^{3}+\left(1-\left\lfloor z^{-1}\right\rfloor z\right)^{3}=8\sigma(x)$,
we get $\left\lfloor z^{-1}\right\rfloor=1$ and $z=\frac{1}{2}+\frac{2}{3}\sqrt{6\sigma(x)-\frac{3}{16}}=\frac12+2\lambda(x)$.
By definition of $g_{\ell}$, we have that for $|x|\leq |b_m|$,
\begin{align}\label{equ:g(sigma)}
g_{\ell}(\sigma(x))=\Big(\frac{1}{4}-\lambda(x)\Big)^{\ell}+\Big(\frac{1}{4}+\lambda(x)\Big)^{\ell}.
\end{align}

Now we apply induction hypothesis to the optimization problem $\OPT_{\ell}^{\sigma'}(m-1,t+1,\rho)$ with variables $r_1,\ldots,r_t,r_{t+1},a_1,b_1,\ldots,a_{m-1},b_{m-1}$, where $r_{t+1}:=a_m$.
Then it follows that $$\rho^{\ell}+\sum_{i=1}^{t}r_i^{\ell}+a_m^{\ell}+\sum_{j=1}^{m-1}\re (a_j+b_j\mi)^{\ell}\geq g_{\ell}(\sigma')=g_{\ell}(\sigma(b_m)).$$
By considering a new function
\begin{align}\label{equ:F(x)}
F_{\ell}(x):=g_{\ell}(\sigma(x))+\re(a_m+x\mi)^\ell-a_m^\ell,
\end{align}
we see from the previous inequality that the objective value of $\OPT_{\ell}^{\sigma}(m,t,\rho)$ satisfies
\begin{align}\label{equ:Phi}
\Phi_{\ell}^{\sigma}(m,t,\rho)=\rho^{\ell}+\sum_{i=1}^{t}r_i^{\ell}+\sum_{j=1}^{m}\re (a_j+b_j\mi)^{\ell}\geq F_{\ell}(b_m).
\end{align}
To prove ($\star$), it further reduces to show $F_{\ell}(b_m)\geq F_{\ell}(0)=g_\ell(\sigma(0))=g_\ell(\sigma)$.	

First, we consider the case when $|b_m|\leq \frac15$.
Note that we have the derivative
\begin{align*}
F'_{\ell}(x)&=\frac{2\ell a_m x}{\lambda(x)}\left[ \left(\frac{1}{4}+\lambda(x)\right)^{\ell-1}-\left(\frac{1}{4}-\lambda(x)\right)^{\ell-1} \right]+\frac{d}{dx}\re(a_m+x\mi)^{\ell}\\
&=\frac{4\ell a_m x}{\lambda(x)}\sum\limits_{j=1}^{2k}\dbinom{\ell-1}{2j-1}\left(\frac{1}{4}\right)^{\ell-2j}\lambda(x)^{2j-1}+\sum\limits_{j=1}^{2k}(-1)^j\dbinom{\ell}{2j}(2j)x^{2j-1}a_m^{\ell-2j},
\end{align*}
where the index $j$ is at most $2k$ as $\ell-1\leq 4k$. Then it follows that for $x\geq 0$,
\begin{align*}
F'_{\ell}(x)&\geq \frac{4\ell a_m x}{\lambda(x)}\sum\limits_{j=1}^{2k-1}\dbinom{\ell-1}{2j-1}\left(\frac{1}{4}\right)^{\ell-2j}\lambda(x)^{2j-1}+\sum\limits_{j=1}^{2k-1}(-1)^j\dbinom{\ell}{2j}(2j)x^{2j-1}a_m^{\ell-2j}\\
&=\ell\cdot\sum\limits_{j=1}^{2k-1}\dbinom{\ell-1}{2j-1}x^{2j-1}a_m^{\ell-2j}\left[\left(\frac{1}{4a_m}\right)^{\ell-2j-1}\left(\frac{\lambda(x)}{x}\right)^{2j-2}+(-1)^j\right].
\end{align*}
Using \eqref{equ:am} we have $0<a_m\leq \frac{1}{80k^2}\leq \frac14$, implying $\left(\frac{1}{4a_m}\right)^{\ell-2j-1}\geq 1$, where each $\ell-2j-1\geq 4k-2-2j\geq 0$.
Also since $\sigma(x)\geq \sigma\geq \frac{73}{800}$, this implies that $\lambda(x)\geq \frac15$ and thus for $|x|\leq |b_m|\leq \frac{1}{5}$, we have $|\frac{\lambda(x)}{x}|\geq 1$.
Putting everything together, we can derive that $F'_{\ell}(x)\geq 0$ for all $|x|\leq |b_m|$.
Therefore, $\Phi_{\ell}^{\sigma}(m,t,\rho)\geq F_\ell(b_m)\geq F_\ell(0)=g_\ell(\sigma)$, as wanted.

\medskip

It remains to consider $|b_m|>\frac{1}{5}$. As a fact we may assume $b_m>\frac15$ (because reversing the sign of $b_m$ would still satisfy all constrictions of $\OPT_{\ell}^\sigma(m,t,\rho)$).
We claim that
\begin{equation}\label{equ:Fbm}
F_{\ell}(b_m)\geq g_{\ell}(\sigma(b_m))-(\ell+1)a_mb_m^{\ell-1}.
\end{equation}
To see this, first consider $\mu\in \{-1,1\}$, that is, $\ell$ is odd. By Lemma~\ref{Re inequality}, we have
\begin{equation*}
F_{\ell}(b_m)\geq g_{\ell}(\sigma(b_m))-a_m^{\ell}-\ell\cdot a_m\left(\sqrt{a_m^2+b_m^2}\right)^{\ell-1}\geq g_{\ell}(\sigma(b_m))-(\ell+1)a_mb_m^{\ell-1},
\end{equation*}
where the last inequality holds because $b_m>\frac{1}{5}, a_m\leq \frac{1}{80k^2}$ and $\ell\geq 5$.
Now consider $\mu=0$ (that is $\ell=4k$).
Let $b_m/a_m=\tan \alpha$ for some $\alpha\in (0,\frac{\pi}{2})$.
Then we have $$\tan \alpha>16k^2>\frac{8k}{\pi}\geq \frac{1}{\tan(\pi/8k)}=\tan\left(\frac{\pi}{2}-\frac{\pi}{8k}\right).$$
So $\alpha\in (\frac{\pi}{2}-\frac{\pi}{8k},\frac{\pi}{2})$ and $\ell\alpha=4k\alpha\in (2k\pi-\frac{\pi}{2},2k\pi)$.
This implies $\re(a_m+b_m\mi)^\ell=(\sqrt{a_m^2+b_m^2})^{\ell}\cos(\ell\alpha)\geq 0$.
Together with \eqref{equ:F(x)}, we also can derive that
$$F_{\ell}(b_m)\geq g_\ell(\sigma(b_m))-a_m^\ell\geq g_{\ell}(\sigma(b_m))-(\ell+1)a_mb_m^{\ell-1},$$
where the last inequality holds since $(b_m/a_m)^{\ell-1}>16k^2>\ell+1$. This proves \eqref{equ:Fbm}.

By \eqref{equ:g(sigma)}, we have
\begin{align*}
\frac{d}{dx}g_{\ell}(\sigma(x))=\frac{2\ell a_m x}{\lambda(x)}\left[ \left(\frac{1}{4}+\lambda(x)\right)^{\ell-1}-\left(\frac{1}{4}-\lambda(x)\right)^{\ell-1} \right]
\end{align*}
For any $x\in [0,b_m]$, we have $\lambda(x)\geq \frac15$ and thus $\left(\frac{1}{4}+\lambda(x)\right)^{\ell-1}\geq \ell\left(\frac{1}{4}-\lambda(x)\right)^{\ell-1}.$
Hence,
\begin{align*}
\frac{d}{dx}g_{\ell}(\sigma(x))\geq \frac{2(\ell-1) a_m x}{\lambda(x)}\left(\frac{1}{4}+\lambda(x)\right)^{\ell-1}.
\end{align*}
Since $w=\frac14+\lambda(x)$ is the larger root of the equation $(\frac{1}{2}-w)^3+w^3=\sigma(x)$, where $\sigma(x)\geq \sigma\geq (\frac{1}{2}-\frac{1}{80k^2})^3+(\frac{1}{80k^2})^3$,
we can deduce that $\frac12>\frac14+\lambda(x)\geq \frac{1}{2}-\frac{1}{80k^2}.$ Then
\begin{align*}
&\frac{d}{dx}g_{\ell}(\sigma(x))\geq 8(\ell-1)a_mx \left(\frac{1}{2}-\frac{1}{80k^2}\right)^{\ell-1}
\end{align*}
for any $x\in [0,b_m]$.	By integration, we have
$g_{\ell}(\sigma(b_m))-g_{\ell}(\sigma)\geq 4(\ell-1)a_mb_m^2\left(\frac{1}{2}-\frac{1}{80k^2}\right)^{\ell-1}.$
This together with \eqref{equ:Phi} and \eqref{equ:Fbm} gives that
\begin{align*}
\Phi_{\ell}^{\sigma}(m,t,\rho)&\geq F_{\ell}(b_m)\geq g_{\ell}(\sigma)+(\ell-1)a_m\left[4b_m^2\left(\frac{1}{2}-\frac{1}{80k^2}\right)^{\ell-1}-\frac{\ell+1}{\ell-1}b_m^{\ell-1}\right]\\
&\geq g_{\ell}(\sigma)+(\ell-1)a_mb_m^{\ell-1}\left[\frac{4}{b_m^{\ell-3}} \left(\frac{39}{80}\right)^{\ell-1}-\frac{\ell+1}{\ell-1}\right]\\
&\geq g_{\ell}(\sigma)+(\ell-1)a_mb_m^{\ell-1}\left[\left(\frac{39}{40}\right)^2\left(\frac{39\sqrt{2}}{40}\right)^{\ell-3}-\frac{\ell+1}{\ell-1}\right]>g_{\ell}(\sigma),
\end{align*}
where the second last inequality follows by that $b_m^2\leq \frac{1}{8}$ and the last inequality holds
because $\left(\frac{39}{40}\right)^2\left(\frac{39\sqrt{2}}{40}\right)^{\ell-3}-\frac{\ell+1}{\ell-1}$ increases for $\ell\geq 5$ and is at least
$2\left(\frac{39}{40}\right)^4-\frac{3}{2}>0$.
This finishes the proof of the statement ($\star$) and thus of Theorem~\ref{4k+3}. \qed

\section{Cycles of length $4k+2$}\label{sec:4k+2}
In this section, we consider cycles $C_\ell$, where $\ell=4k+2$ for some integer $k\geq 1$.

First, we construct a family of tournaments $T$ with $t(C_\ell, T)< g_\ell(t(C_3,T))$ for any value of $t(C_3,T)$.
We need to introduce some of the limit theory of tournaments established in \cite{Ckgood}.
A {\it tournamenton} is a measurable function $W: [0,1]^2\to [0,1]$ such that $W(x,y)+W(y,x)=1$ for all $(x,y)\in [0,1]^2$.
For a tournamenton $W$, define $$C(W,\ell)=2^\ell \int_{x_1,...,x_\ell\in [0,1]} W(x_1,x_2)W(x_2,x_3)\dots W(x_{\ell-1},x_\ell)W(x_\ell,x_1)~dx_1\dots x_\ell.$$
One can also define the {\it spectrum} $\sigma(W)$ (for its precise definition we refer to Section 2.2 of \cite{Ckgood}).
By Proposition~4 of \cite{Ckgood}, $C(W,\ell)=2^\ell\cdot \sum_{x\in \sigma(W)} x^\ell.$
A {\it carousel tournament} $T_n$ is a tournament with vertex set $\{0,1,...,2n\}$, where $i\to j$ for every $i, j$ satisfying $i+1\leq j\leq i+n$ (computations module $2n+1$).
Note that $T_n$ is regular, implying that $t(C_3,T_n)=1/8+o(1)$ where $o(1)\to 0$ as $n\to \infty$.
Let $W_C$ be the tournamenton as follows: for $x,y\in [0,1]$,
let $W_C(x,x)=1/2, W_C(x,y)=1$ if $y\in (x-1,x-1/2)\cup (x,x+1/2]$, and $W_C(x,y)=0$ otherwise.
Then $W_C$ is the limit object of carousel tournaments $T_n$ satisfying that $\lim_{n\to \infty}t(C_\ell,T_n)=C(W_C,\ell)/2^\ell$.
One can deduce from \cite{Ckgood} that $$\sigma(W_C)=\{1/2,0\}\cup \{\pm \mi/(2k-1)\pi: \forall k\in \mathbb{N}\}$$
and as $\ell=4k+2$,
\begin{align}\label{equ:t(4k+2)}
\lim_{n\to \infty}t(C_\ell,T_n)=\frac{C(W_C,\ell)}{2^\ell}=\sum_{x\in \sigma(W_C)}x^\ell=\frac{1}{2^\ell}-2\sum_{k=1}^\infty \Big(\frac{1}{(2k-1)\pi}\Big)^\ell\triangleq\alpha_\ell.
\end{align}
Note that $\alpha_\ell<1/2^\ell$. We are ready to construct the desired family of tournaments.

\begin{defi}
A {\it carousel blow-up} of a $m$-vertex transitive tournament is a tournament $T$ with $V(T)=V_1\cup V_2\cup ...\cup V_m$
such that each $V_i$ induces a carousel tournament and for any $i<j$, all arcs between $V_i$ and $V_j$ are oriented from $V_i$ to $V_j$.
A carousel blow-up of a $m$-vertex transitive tournament is called {\it balanced} if $|V_1|=|V_2|=\dots =|V_{m-1}|\geq |V_m|$.
\end{defi}

Now consider a balanced carousel blow-up $T^\star$ of a transitive tournament, and compute the density $t(C_\ell, T^\star)$.
Suppose that $T^\star$ has $n$ vertices and $t$ parts of equal size $zn$ for some $z\in (0,1]$.
Then $tz\leq 1< (t+1)z$, which implies that $t=\lfloor z^{-1}\rfloor$.
By \eqref{equ:t(4k+2)}, if $n$ is sufficiently large, then
$$t(C_\ell, T^\star)=\alpha_\ell\cdot \Big(\lfloor z^{-1}\rfloor z^\ell+(1-\lfloor z^{-1}\rfloor z)^\ell\Big)+o(1), \mbox{ ~~ where } o(1)\to 0 \mbox { as } n\to \infty.$$
On the other hand, we have $t(C_3,T^\star)=\frac18\Big(\lfloor z^{-1}\rfloor z^3+(1-\lfloor z^{-1}\rfloor z)^3\Big)+o(1)$.
Note that $z$ can be any real in $(0,1]$.
We can summarize this construction as following.

\begin{lem}\label{lem:4k+2}
For any $\ell=4k+2$ for some integer $k\geq 1$, there exist tournaments $T$ with arbitrary $t(C_3,T)$ such that
$t(C_\ell, T)=2^\ell \alpha_\ell \cdot g_\ell(t(C_3,T))+o(1)$, which is strictly less than $g_\ell(t(C_3,T)).$
\end{lem}

Next, we apply results of \cite{Ckgood} to obtain a lower bound on $t(C_\ell,T)$ for all regular tournaments $T$.
We need the following lemma in \cite{Ckgood}.
Let $D_n$ denote the $n\times n$ skew-symmetric matrix with all entries above the diagonal equal to $1$ and all entries below the diagonal equal to $-1$.
\begin{lem}[\cite{Ckgood}, Lemma~11]\label{Dn}
For any $n\in \mathbb{N}$, the spectral radius of any skew-symmetric matrix in $[-1,1]^{n\times n}$ is at most the spectral radius of $D_n$.
\end{lem}

Let $A$ be the tournament matrix of an $n$-vertex regular tournament $T$.
Let $J$ be the $n\times n$ matrix with every entry equal to $1$ and let $B:=A-J/2n$.
Then $B$ is a skew-symmetric matrix with all entries in $\{\frac{-1}{2n},\frac{1}{2n},0 \}$.
Because $T$ is regular, the sum of entries in each column (or row) of $B$ is $0$, which shows that $JB=BJ=0$.
Therefore by \eqref{equ:trA}, we have $$t(C^\ell,T)=tr(A^{\ell})+o(1)=tr\left((J/2n+B)^\ell\right)+o(1)=tr((J/2n)^\ell)+tr(B^\ell)+o(1).$$
It is known that the spectral radiuses $\rho_n$ of the matrices $D_n$ divided by $n$ converge to $\frac{2}{\pi}$ (see the proof of Lemma~12 in \cite{Ckgood}).
Let the spectral radius of $B$ be $\rho$ and  the eigenvalues of $B$ be $\lambda_1,\lambda_2,\ldots,\lambda_n$.
Then we have that by Lemma~\ref{Dn}, $|\lambda_i|\leq \rho\leq \frac{\rho_n}{2n}=\frac{1}{\pi}+o(1)$ for any $i\in [n]$,
and $\sigma\triangleq\sum\limits_{i=1}^{n}|\lambda_i|^2\leq tr(BB^T)=\frac{n-1}{4n}$,
where the last equality follows from the fact that every non-diagonal entry of $B$ is $1/2n$ or $-1/2n$ and every diagonal entry of $B$ is 0.
Applying Lemma~\ref{n}, as $\ell$ is even, we can obtain that $$|tr(B^\ell)|\leq \sum\limits_{i=1}^{n}|\lambda_i|^{\ell}=\left\lfloor\frac{\sigma}{\rho^2}\right\rfloor\rho^{\ell}+\left(\sigma-\left\lfloor\frac{\sigma}{\rho^2}\right\rfloor\rho^2\right)^{\ell/2}.$$
Note that $\lfloor\frac{\sigma}{\rho^2}\rfloor\rho^{\ell}+(\sigma-\lfloor\frac{\sigma}{\rho^2}\rfloor\rho^2)^{\ell/2}$ increases as $\sigma$ and $\rho$ increase.
By plugging $\sigma\leq \frac{n-1}{4n}$ and $\rho\leq\frac{1}{\pi}+o(1)$ in the above inequality,
we obtain $$|tr(B^\ell)|\leq 2\left(\frac{1}{\pi}\right)^\ell+\left(\frac{1}{4}-\frac{2}{\pi^2}\right)^{\ell/2}+o(1)$$
and thus $$t(C^\ell,T)\geq tr\left((J/2n)^\ell\right)-|tr(B^\ell)|+o(1)\geq \frac{1}{2^\ell}- 2\left(\frac{1}{\pi}\right)^\ell-\left(\frac{1}{4}-\frac{2}{\pi^2}\right)^{\ell/2}+o(1),$$
where $\left(\frac{1}{4}-\frac{2}{\pi^2}\right)^{1/2}\approx 0.218 <\frac1{\pi}$.
Together with \eqref{equ:t(4k+2)}, we have the following lemma.

\begin{lem}\label{lem:regular-4k+2}
For every $\epsilon>0$, there exists $\ell_0$ such that for every $\ell\geq \ell_0$ with $\ell\equiv 2\mod 4$,
$$\frac{1}{2^\ell}-(2+\epsilon)\left(\frac{1}{\pi}\right)^\ell \leq \min_{T} t(C_\ell, T)\leq  \frac{1}{2^\ell}-2\left(\frac{1}{\pi}\right)^\ell,$$
where the minimum is over all $n$-vertex regular tournaments $T$ for large $n$.
\end{lem}

It seems plausible to believe that for every such $\ell\geq 6$, $\lim_{n\to \infty}\min_{T} t(C_\ell, T)=\alpha_\ell$ where the minimum is over all $n$-vertex regular tournaments $T$.

\section{Concluding remarks}
We prove the statement of Theorem~\ref{4k+3} whenever $t(C_3,T)\geq \frac18-O(\frac{1}{\ell^2})$.
It is possible to lower the condition on $t(C_3,T)$.
A special case is $\ell=5$.
Using Theorem~\ref{question} and a similar argument as in \cite{C3C4} (via the method of Lagrange multiplier),
we also can obtain the same bound as in Theorem~\ref{C3C4theo} that $t(C_5,T)\geq g_5(t(C_3,T))+o(1)$ whenever $t(C_3,T)\in [\frac{1}{72},\frac{1}{8}]$.
For integers $\ell$ with $\ell\equiv 0,1 \mod 4$, using more careful calculation,
we can prove the same conclusion of Theorem~\ref{4k+3} whenever $t(C_3,T)\geq \frac18-O(\frac{1}{\ell})$.
However, the approach we used has its obvious limit, so we choose to present a unified, less complicated proof for Theorem~\ref{4k+3}.

We now conclude this paper by the following two conjectures.
The first one is a direct generalization of Conjecture~\ref{conj1} for cycles $C_\ell$ with $\ell\not\equiv 2\mod 4$,
which is supported by Theorem~\ref{4k+3}.

\begin{conj}
For any integer $\ell\geq 4$ with $\ell\not\equiv 2\mod 4$,
every tournament $T$ satisfies that
$$t(C_\ell,T)\geq g_\ell(t(C_3,T))+o(1),$$
where the $o(1)$ term goes to zero as $|V(T)|$ goes to infinity.
\end{conj}

\begin{conj}
For any integer $\ell\geq 6$ with $\ell\equiv 2\mod 4$,
every tournament $T$ satisfies that
$$t(C_\ell,T)\geq 2^\ell \alpha_\ell \cdot g_\ell(t(C_3,T))+o(1),$$
where the $o(1)$ term goes to zero as $|V(T)|$ goes to infinity.
\end{conj}

\noindent If true, as demonstrated by Lemma~\ref{lem:4k+2}, the above conjecture would be sharp.

\medskip

\bigskip

\noindent {\bf Acknowledgement.}
The authors would like to thank Prof. Kimchuan Toh for many helpful discussions and valuable comments which greatly improve the presentation.

\end{document}